\documentstyle[11pt]{article}
\begin{document}

\newtheorem{lemma}{Lemma}[section]
\newtheorem{theorem}[lemma]{Theorem}
\newtheorem{definition}[lemma]{Definition}
\newtheorem{observation}[lemma]{Observation}
\newtheorem{remark}[lemma]{Remark}
\newtheorem{corollary}[lemma]{Corollary}
\newtheorem{proposition}[lemma]{Proposition}
\newtheorem{example}[lemma]{Example}
\newtheorem{algorithm}[lemma]{Algorithm}
\newtheorem{conjecture}[lemma]{Conjecture}

\title{\bf On the Complexity of \\
Smooth Projective Toric Varieties}
\author{
 Serkan Ho\c{s}ten \\
School of OR \& IE, Cornell University, Ithaca, NY, 14853\\
{\tt serkan@cs.cornell.edu}\\
\phantom{dadada} 
} 
\maketitle

\begin{abstract}
In this paper we answer a question posed by V.V.~Batyrev in
\cite{Bat}. The question asked if there exists a complete
regular fan with more than {\bf quadratically} many primitive collections.   
We construct a smooth projective toric variety associated to a complete 
regular fan $\Delta$ in ${\bf R}^d$ with $n$ generators where the number
of primitive collections of $\Delta$ is at least {\bf exponential} in $n-d$. 
We also exhibit the connection between the number of 
primitive collections of $\Delta$ and the facet complexity
of the Gr\"obner fan of the associated integer program.  
\end{abstract}

\vskip 0.7cm
\section{Introduction}

\noindent
In this paper we give an affirmative answer to the following question
posed by V.V.Batyrev \cite{Bat}:

\vskip 0.3cm
\noindent
{\bf Question:} {\sl
Does there exist a complete regular $d$-dimensional fan $\Delta$ 
with $n$ generators such that $\Delta$ has more than
$(n-d-1)(n-d+2)\slash2$ primitive collections for $n-d > 1$?
}

\vskip 0.3cm
\noindent
In Section 2 we prove the following theorem which answers the above
question.
\begin{theorem}
There exists a complete regular $d$-dimensional fan $\Delta$ with $n$
generators where the number of primitive collections of $\Delta$ is 
more than $2^{\frac{1}{2}(n-d)}$.
\end{theorem}

\vskip 0.3cm
\noindent
A fan $\Delta \subset {\bf R}^d$ that covers ${\bf R}^d$ is a complete
fan. If we require the full-dimensional cones in $\Delta$ to
be simplicial with integral generators which form a ${\bf Z}$-basis
for ${\bf Z}^d$, then $\Delta$ is said to be {\it regular} (see below
for formal definitions).   
If a fan $\Delta \subset {\bf R}^d$ is generated by the
$1$-dimensional cones defined by the vectors in ${\cal A}
=  \{a_1, a_2, \ldots, a_n\} \subset {\bf Z}^d$, 
then the {\it primitive collections} of
$\Delta$ are defined as follows:
\begin{definition}
A nonempty subset ${\cal P} = \{a_{i_1},a_{i_2},\ldots,a_{i_k}\}$ of 
${\cal A}$ is called a
{\bf primitive collection} if for each generator $a_{i_s} \in {\cal P}$
the elements ${\cal P} \setminus a_{i_s}$ generate a $(k-1)$-dimensional
cone in $\Delta$, while ${\cal P}$ does not generate a $k$-dimensional
cone in $\Delta$.
\end{definition}

\noindent
The primitive collections of a complete regular 
fan $\Delta \subset {\bf R}^d$ with $n$ generators are studied in
\cite{Bat} to classify $d$-dimensional smooth complete toric varieties
with $n-d=3$.
The question we study asks whether there exists a complete regular 
fan $\Delta$ where the number of primitive collections of $\Delta$ is
at least {\it quadratic} in $n-d$. Theorem 1.1 constructs a complete regular
fan with {\it exponentially} many primitive collections. 
The same theorem can be restated in the language of Gr\"obner bases
of toric varieties: 
Theorem 3.1 shows that there exists a toric variety  $X$
with a square-free initial ideal whose number of
minimal generators is exponential in the codimension of $X$.
In the last section
we make a connection 
between two conjectures: one of them appears in the context of 
the complexity of complete regular fans (Conjecture 7.1 in \cite{Bat})
and the other one is about the complexity of {\it Gr\"obner fans} in
the context of {\it integer programming} (Conjecture 6.1 in
\cite{ST}).  

\noindent
In this article we will use results from the theory of 
{\it coherent triangulations} of a vector configuration. In order to prove
Theorem 1.1
we need the following definitions which connect
coherent triangulations
and complete regular projective fans: 
A complete fan $\Delta \in {\bf R}^d$
with $n$ generators ${\cal A}
= \{a_1, a_2, \ldots, a_n\} \subset {\bf Z}^d$ is said to be {\it regular}
if every $d$-dimensional cone $\sigma \in \Delta$ is simplicial and 
the generators $\{a_{i_1}, a_{i_2}, \ldots, a_{i_d}\}$ of $\sigma$
form a ${\bf Z}$-basis of ${\bf Z}^d$.

\begin{definition}
A complete regular fan $\Delta$ is said to be {\bf projective} if there exists
a {\bf support function} $\phi : {\bf R}^d \rightarrow {\bf R}$ such that
\begin{enumerate}
\item $\phi$ is convex and $\phi({\bf Z}^d) \subset {\bf Z}$,
\item $\phi$ is linear on each cone of $\Delta$ with $\phi\mid_{\sigma}
\neq \phi\mid_{\tau}$ for distinct $d$-dimensional cones $\sigma$
and $\tau$.
\end{enumerate}
\end{definition}

\noindent
It is a well-known fact that if $V(\Delta)$ is the smooth complete 
$d$-dimensional toric variety that is associated with $\Delta$ then 
$V(\Delta)$ is projective if and only if $\Delta$ is projective
\cite{Oda}. In this paper we will use an equivalent definition of 
a projective toric variety via {\it coherent} triangulations. 

\begin{definition} A triangulation $T$ of a vector configuration ${\cal A}
= \{a_1, a_2, \ldots, a_n\} \in {\bf R}^d$ is a polyhedral complex 
consisting of simplical cones which cover 
$pos({\cal A}) = \{x \in {\bf R}^d\,:
\, x = \sum_{i=1}^{n} \lambda_i a_i, \, \lambda_i \geq 0\}$.
A triangulation $T$ of ${\cal A}$ is said to be {\bf coherent} 
if there exists a support function $\phi$ on $T$ as in Definition 1.3 (see
\cite{BFS}, \cite{GKZ}).
\end{definition}
 
\section{Exponential lower bound}

\noindent
To give an exponential lower bound for the number of primitive
collections of a complete regular fan we will use the example
given in Proposition 6.7 of \cite{ST}. 
\begin{definition}
Given a matrix $B \in {\bf Z}^{d \times n}$, the {\bf chamber complex} 
$\Gamma(B)$ of $B$ is the coarsest polyhedral complex that refines
all triangulations of $B$ and covers $pos(B)$.
\end{definition}
\begin{proposition}
\cite{BGS} Let $A$ be a Gale transform of 
$B \in {\bf Z}^{d \times n}$, i.e. let 
$A \in {\bf Z}^{(n-d) \times n}$ such that
$$0 \longrightarrow {\bf R}^{d} \stackrel{B^T}{\longrightarrow}
 {\bf R}^{n} \stackrel{A}{\longrightarrow} {\bf R}^{(n-d)} 
\longrightarrow 0$$is exact. Then there is a bijection between the 
coherent triangulations of $A$ and the $d$-dimensional chambers of $\Gamma(B)$
given by
$$T = \bigcup_{i=1}^{t} \sigma_i  \Longleftrightarrow  C =
\bigcap_{i=1}^{t} \sigma^{\ast}_i$$
where $\sigma_i = pos(a_{i_1},a_{i_2},\ldots,a_{i_{n-d}})$ are the cones of
the coherent triangulation $T$ and $\sigma^{\ast}_i = pos(\{b_j: 1 \leq
j \leq n, \,j \neq i_1, i_2, \ldots, i_{n-d}\})$ are the cones of $B$ 
containing the chamber $C$.
\end{proposition}

\noindent Now we construct a complete regular projective fan which
has exponentially many primitive collections.
Let $B$ be the node-edge incidence matrix of the complete
bipartite graph $K_{n,m}$ where $n=2k-1$ and $m=2k+1$.
$B = \{e_i \times e'_j: \, 1 \leq i \leq n, 1 \leq j \leq m \}$ 
where $e_i \in {\bf R}^n$ and $e'_j \in {\bf R}^m$ are unit vectors.
$B$ has rank $n+m-1$ and is unimodular, i.e. any subdeterminant of
$B$ is $0$ or $\pm 1$ (\cite{Schr}, p.273).
The cone $pos(B)$  consists of all non-negative
vectors $(u_1,\ldots,u_n) \times (v_1,\ldots,v_m)$ such that
$u_1+\cdots+u_n = v_1+\cdots+v_m$.
Let $A \in {\bf Z}^{(n-1)(m-1) \times nm}$ be a
Gale transform of $B$. By Proposition 2.2, for every chamber in 
the chamber complex $\Gamma(B)$ there exists a corresponding 
coherent triangulation of $A$. We consider a special chamber in
$\Gamma(B)$. The one-dimensional cone
generated by $(\frac{1}{n},\ldots,\frac{1}{n}) \times 
(\frac{1}{m},\ldots,\frac{1}{m})$ is in $pos(B)$ and we claim that
it is in the interior of a full-dimensional chamber.
The facets of full-dimensional chambers of $\Gamma(B)$ 
correspond to the cocircuits of the oriented matroid of $B$ 
(\cite{DHSS}, Lemma 2.7). In our situation, a cocircuit 
of $B$ corresponds to a cut $(C_+,C_-;D_+,D_-)$ in $K_{n,m}$ where
$(C_+,C_-)$ is a partition of $\{1,\ldots,n\}$ and 
$(D_+,D_-)$ is a partition of $\{1,\ldots,m\}$. The corresponding
hyperplane is defined by
$$\sum_{i \in C_+} u_i - \sum_{i \in C_-} u_i - \sum_{j \in D_+} v_j +
\sum_{j \in D_-} v_j \, = \, 0. \eqno (1)$$
Since $n$ and $m$ are relatively prime,
$(\frac{1}{n},\ldots,\frac{1}{n}) \times
(\frac{1}{m},\ldots,\frac{1}{m})$ cannot lie on any of these
hyperplanes. So it must be in the interior of a full-dimensional
chamber. This chamber is called the {\it central chamber} of $pos(B)$.
Let $\Delta$ be the corresponding coherent triangulation of $A$.
Since $B$ is unimodular, so is $A$, and therefore 
$\Delta$ consists of 
simplical cones whose generators form a ${\bf Z}$-basis for 
${\bf Z}^{(n-1)(m-1)}$. Because there exists a strictly positive
vector in $im(B^T) = ker(A)$, any vector in ${\bf
R}^{(n-1)(m-1)}$ is in $pos(A)$. This shows that $\Delta$ is 
a complete regular fan. 

\noindent
Now we show that every column of $A$ is a generator of
$\Delta$. By Proposition 2.2 it is enough to show that for every
column $b_{i,j} = e_i \times e'_j \in B$ there exists 
a cone $\tau$ which contains 
the central chamber but which does not have $b_{i,j}$ as a
generator. Suppose $\tau$ is a cone that contains the central chamber
and has $b_{i,j}$ as a generator. Since the natural action
of the product of symmetric groups $S_n \times S_m$ on $B$ fixes the
central chamber, for any $\pi \times \sigma \in S_n \times S_m$, 
$(\pi \times \sigma)(\tau)$ covers the central chamber as well. We can
pick $\pi$ such that $\pi(i) = i$. If $b_{i,k}, \, k \neq j$ does not 
appear as a generator of $\tau$ we can choose $\sigma$ such that
$\sigma(k) = j$, and we would be done. So suppose $b_{i,k}$ is a
generator of $\tau$ for $k = 1,\ldots,m$. Since these vectors are
linearly independent and $rank(B)=n+m-1$, there are exactly $n-1$ 
generators of $\tau$ which are not of the above form. But 
$\tau$ contains $(\frac{1}{n},\ldots,\frac{1}{n}) \times
(\frac{1}{m},\ldots,\frac{1}{m})$ in its interior, so these remaining
$n-1$ generators are of the form $b_{k,s_k}, \, k =1,\ldots,n, \, k \neq
i$. Each of these generators must have the coefficient $\frac{1}{n}$
in the unique expression that expresses 
$(\frac{1}{n},\ldots,\frac{1}{n}) \times (\frac{1}{m},\ldots,\frac{1}{m})$
in terms of generators of $\tau$. But $\frac{1}{n} > \frac{1}{m}$, and
this shows that $\Delta$ is a complete regular projective fan
generated by the columns of $A$. 

\noindent
In order to give an exponential lower bound on the primitive collections of
$\Delta$ constructed above we establish 
a link between its {\it circuits} and its primitive collections. 
\begin{definition} Let ${\cal A}$ be a vector configuration in ${\bf
Z}^d$. A collection of linearly dependent vectors 
$Z \subseteq {\cal A}$ is called a {\bf circuit} 
if any proper subset of $Z$ is linearly independent. 
\end{definition}
   
\noindent
We will call the circuits of the generators of a complete fan $\Delta$
the circuits of $\Delta$. If $Z$ is a circuit of $\Delta$, the unique
(up to sign) dependence relation $ \sum_i \lambda_i z_i = 0$
partitions $Z$ into two subsets, namely $Z_+ = \{z_i \in Z:
\lambda_i > 0\}$ and $Z_i = \{z_i \in Z:\lambda_i < 0\}$. In this
case, there exist
precisely two triangulations of $Z$: $\, t_+(Z) = \{Z \setminus
z_i: z_i \in Z_+\}$ and $t_-(Z) = \{Z \setminus z_i: z_i \in Z_-\}$.
Note that $relint(pos(Z_+)) \bigcap relint(pos(Z_-)) \neq \emptyset$
(see \cite{BLSWZ}).
Given a triangulation $\Delta$ and a circuit $Z$ of $\Delta$ such that
$t_+(Z)$ is a subcomplex of $\Delta$, one can get via a  
{\it bistellar flip} another triangulation
$\Delta'$ such that $t_-(Z)$ is a subcomplex of $\Delta'$. For the
details we refer to [GKZ, pp. 231-233]. 
The next lemma makes the connection
between the circuits and primitive collections of $\Delta$.
\begin{lemma}
Let $\Delta \subset {\bf R}^d$ be a complete regular fan (i.e. a
 triangulation) and let $Z$ be a circuit
such that $t_+(Z)$ is a subcomplex of $\Delta$. Then $Z_+$ is a
primitive collection. Moreover, if $Z'$ is a different
circuit where $t_+(Z')$ is a subcomplex of $\Delta$, then $Z_+ \neq
Z'_+$.
\end{lemma}

\noindent {\sl Proof:} Clearly $Z_+$ does not generate a cone 
in $\Delta$. By definition of $t_+(Z)$, for all $z \in Z_+, \, pos(Z_+
\setminus z)$ is a face of $t_+(Z)$, and hence is a cone in $\Delta$. 
Each of these cones must be $(card(Z_+)-1)$-dimensional, since otherwise
$Z$ cannot be a circuit.
This shows $Z_+$ is a primitive collection. For the second statement,
assume $Z_+ = Z'_+$. Since $t_+(Z) \neq t_+(Z')$, the respective 
subcomplexes $K$ and $K'$ of $\Delta$ 
on which the bistellar flips are supported are different as well. 
But $pos(Z_+) \bigcap relint(K) \neq \emptyset$ and 
$pos(Z'_+) \bigcap relint(K') \neq \emptyset$, which implies
$relint(K) \bigcap relint(K') \neq \emptyset$. This cannot happen 
since $K$ and $K'$ are distinct subcomplexes of $\Delta$. This
contradiction completes the proof. $\Box$

\vfill
\eject
\noindent 
For the main theorem we need the following result 
which provides the link between bistellar flips (and hence the
primitive collections) of $\Delta$ and the corresponding 
chamber in the dual configuration.

\begin{theorem}(\cite{GKZ}, p.233) Let ${\cal A} \subset {\bf Z}^d$ be a 
vector configuration and let ${\cal B}$ be a Gale transform of ${\cal
A}$. If $\Delta$ and $\Delta'$ are two coherent triangulations of ${\cal A}$,
then $\Delta$ and $\Delta'$ differ by a bistellar flip if and only if 
the corresponding chambers in $\Gamma({\cal B})$ share a facet.
\end{theorem}

\noindent {\sl Proof of Theorem 1.1:}
Let $B$ be the node-incidence matrix of $K_{n,m}$ where $n=2k-1$ and
$m=2k+1$ and let $A$ be a Gale transform of $B$. Let $\Delta$ be the
coherent triangulation of $A$ that corresponds to the central chamber
in $\Gamma(B)$. As we established before $\Delta$ is a complete regular 
projective fan generated by the columns of $A$. Lemma 2.4 and Theorem
2.5 imply that the number of
primitive collections of $\Delta$ should be at least the number of
facets of the central chamber in $\Gamma(B)$. 
The following proposition shows that there are at least exponentially
many such facets. The proof of the proposition 
can be found in \cite{ST}, but we include its proof for completeness.

\begin{proposition}
The central chamber in $\Gamma(B)$ which correponds to $\Delta$ has 
at least $4^k$ facets.
\end{proposition}

\noindent {\sl Proof:} If $H$ is a hyperplane defined by the 
equation (1), we will call $(card(C_+),card(D_+))$ the {\it type} of
$H$. Now starting at the point $(\frac{1}{n},\ldots,\frac{1}{n}) \times
(\frac{1}{m},\ldots,\frac{1}{m})$ and moving in the direction 
of $(-1,\ldots,-1,n-1) \times (0,\ldots,0)$ to a generic 
point $(a,\ldots,a,a+1-na) \times (\frac{1}{m},\ldots,\frac{1}{m})$
we cross a facet of type $(r,s)$ with $n \in C_-$ whenever
$$\, r \cdot a \,\, - \,\, (n-r-1) \cdot a - (a+1-na) \,\, - \,\, 
s / m \,\, + \,\, (m-s)/m \quad = \quad 
2 \cdot r \cdot a \,-\, 2 \cdot s / m 
\quad = \quad 0. $$
From here we get $a = s /  mr$ and since $a < 1 / n$ we like 
to find $r$ and $s$ which minimize the positive integer $m \cdot r - n
\cdot s$. The unique solution is $r = k$ and $s = k + 1$ and since 
$S_n \times S_m$ acts transitively on the set of hyperplanes of 
type $(r,s)$, we conclude that every hyperplane of this type is a 
facet of the central chamber. When $k \geq 3$ there are $2k-1 \choose k$
$2k+1 \choose k+1$ $> 2^k \cdot 2^k = 4^k$ such facets. $\Box$
 
\section{Connections to integer programming}

\noindent
In this section we will first state Theorem 1.1 in terms of the {\it Gr\"obner 
basis} (\cite{AL}, \cite{CLO}) of an {\it integer
program}. Subsequently we will relate two conjectures, one that
appears in the context of smooth projective toric varieties and the
other one in the context of integer programming. 
An integer
program can be stated as follows:
$$ \hbox{minimize} \,\, c \cdot x \,\,\,\,
\hbox{subject to} \,\, A \cdot x \, = \, b,
 \,\,  x \in {\bf N}^n$$
where $A \in {\bf Z}^{d \times n}$ with $rank(A) = d$, $b \in {\bf
Z}^d$ and $c \in {\bf R}^n$. The reduced Gr\"obner basis of the 
{\it toric ideal} $I_A = \langle x^\alpha - x^\beta: \, \alpha,
\beta \in {\bf N}^n, \, A \cdot \alpha = A \cdot \beta \rangle$
with respect to the term order induced by the cost vector $c$ provides
a {\it test set} for solving this integer program (see \cite{AL},
\cite{CT}, \cite{ST}, \cite{Th} for details).
\begin{theorem}
Let $A$ be a Gale transform of the node-edge incidence matrix $B$ of 
$K_{n,m}$ with $n=2k-1$ and $m=2k+1$. Then 
$I_A$ is generated by $x_{i1} x_{i2} \cdots x_{im} - 1, \, i = 1,\ldots,n$  and
$x_{1j} x_{2j} \cdots x_{nj} -~1, \, j = 1, \ldots,m$, 
and the reduced Gr\"obner
basis of $I_A$ with respect to the degree lexicographic term order contains
at least $4^k$ elements.
\end{theorem}

\noindent {\sl Proof:} 
The rows of $B$ constitute a ${\bf Z}$-basis for $ker(A) \bigcap {\bf
Z}^{nm}$. The above binomials correspond to the rows of $B$ and the
ideal they generate is contained in $I_A$. But since the sum of all
the rows of $B$ is a strictly positive vector, these binomials
generate $I_A$ (see Lemma 2.1 in \cite{SWZ}).
The degree lexicographic term order $\succ_{deglex}$ can be represented by the 
cost vector $c = (1,1,\ldots,1)$ refined by the lexicographic order.
Since $A$ is unimodular, the reduced Gr\"obner basis of $I_A$ with
respect to $\succ_{deglex}$ consists of 
square-free binomials (\cite{St}, Corollary 8.9) and the initial 
term of each binomial corresponds to a minimal non-face (i.e. a
primitive collection) of
the coherent triangulation $\Delta$ induced by $c$ (\cite{St}, Theorem
8.3). As the vector $B\cdot c$ is in the central chamber of
$\Gamma(B)$, Proposition 2.2 implies that  the triangulation $\Delta$
induced  by $c$ is the
same as the complete regular fan we considered in the proof of Theorem
1.1. This finishes the proof. $\Box$ 

\begin{example}($3 \times 5$ Complete Bipartite Graph)
\end{example} Let $B$ be the node-incidence matrix of $K_{3,5}$. If we
associate with every column $b_{ij} = e_i \times e^{\prime}_j, \, \,
i = 1,2,3, \,\, j = 1,2,3,4,5$, the variable $x_{ij}$, 
then 
$$I_A \, = \, \langle x_{11}x_{21}x_{31} - 1, x_{12}x_{22}x_{32} - 1, 
x_{13}x_{23}x_{33} - 1, x_{14}x_{24}x_{34} - 1, x_{15}x_{25}x_{35} - 1,$$
$$x_{11}x_{12}x_{13}x_{14}x_{15} - 1, 
x_{21}x_{22}x_{23}x_{24}x_{25} - 1, x_{31}x_{32}x_{33}x_{34}x_{35} - 1
\rangle.$$
There are 30 facets of the central chamber of $B$ and indeed the reduced
Gr\"obner basis of $I_A$ with respect to $\succ_{deglex}$ consists of
50 binomials.

\vskip 0.3 cm 
\noindent
In relation to smooth complete projective
varieties, the following conjecture is posed in \cite{Bat}.
\begin{conjecture} (\cite{Bat}, Conjecture 7.1)
For any $d$-dimensional smooth complete toric variety  defined by a
complete regular fan $\Delta$ with $n$ generators, there exists a constant
$N(n-d)$ depending only on $n-d$ such that the number of primitive 
collections in $\Delta$ does not exceed $N(n-d)$.
\end{conjecture}
Another conjecture with a similar flavor is stated about the
complexity of {\it Gr\"obner cones} in the setting of integer
programming in \cite{ST}
(Conjecture 6.1), and here we will give a connection between the two
conjectures along the lines of the previous section. Given an integer
program defined by a matrix $A$,
two generic cost vectors 
$c$ and $c'$ are considered to be equivalent if the respective reduced
Gr\"obner bases of $I_A$ are the same. The set of all such 
equivalent cost vectors
associated to a fixed reduced Gr\"obner basis of the toric ideal of
$A$ is an open polyhedral cone and the collection of the closures 
of all such cones and their faces
constitute a fan called the {\it Gr\"obner fan} of $A$ (\cite{MR}, \cite{BM},
\cite{St}, \cite{ST}). 
\begin{conjecture} (\cite{ST}, Conjecture 6.1)
There exists a function $\varphi$ such that, for every matrix $A \in 
{\bf Z}^{d \times n}$ of rank $d$, every cone of the Gr\"obner fan of 
$A$ has at most $\varphi(n-d)$ facets.
\end{conjecture}
This conjecture is true for $n-d \leq 2$. For the case $n-d = 3$, 
$\varphi(3) = 4$ under certain genericity assumptions on the matrix
$A$ and this was proved by I.~B\'ar\'any and H.~Scarf in \cite{BS}.
The following proposition points to a connection between $N$ and
$\varphi$:
\begin{proposition} (\cite{ST}, Corollary 3.18) Let $A \in 
{\bf Z}^{d \times n}$ with $rank(A)=d$ be a unimodular matrix and 
let $B$ be a Gale transform of $A$. Then the Gr\"obner fan of $A$ and
$\Gamma(B)$ coincide.
\end{proposition}
In the light of this proposition one can formulate a specialized version
of Conjecture 3.4.
\begin{conjecture} 
There exists a function $\varphi'$ such that, for every unimodular
matrix 
$A \in {\bf Z}^{d \times n}$ of rank $d$, every cone of $\Gamma(B)$  
has at most $\varphi'(n-d)$ facets where $B$ is a Gale transform of $A$.
\end{conjecture}
\begin{theorem} If there exist $N$ and $\varphi'$ as above, then 
$N(n-d) \geq \varphi'(n-d)$ for all $n$ and $d$.
\end{theorem}

\noindent {\sl Proof:} Fix $n$ and $d$, and suppose that $A \in 
{\bf Z}^{d \times n}$ is a unimodular matrix with a coherent
triangulation $\Delta$ such that the corresponding chamber in
$\Gamma(B)$, where $B$ is a Gale transform of $A$, has $\varphi'(n-d)$
facets. If $\Delta$ uses all columns of $A$ as generators, then by the
results of the previous section we would be done. Otherwise we can 
refine $\Delta$ into another complete regular fan by adding the 
missing generators. This will not destroy the primitive collections in
$\Delta$ associated with the bistellar flips as in Lemma 2.4. Hence
$N(n-d) \geq \varphi'(n-d)$. $\Box$  

\vskip 0.5cm
\noindent {\bf Acknowledgement} The author thanks Bernd
Sturmfels for pointing out Batyrev's open problem.

\end{document}